\documentclass[12pt]{article}
\usepackage{amssymb}
\date{}
\makeatletter

\@addtoreset{equation}{section}
\makeatother
\newcommand{\qed}{\hbox{\rule[-2pt]{3pt}{6pt}}}

\begin{document}
\title {\bf Asymptotic behavior of bifurcation curves of ODEs with oscillatory nonlinear diffusion 
 }

\author{Tetsutaro Shibata
\\
Laboratory of Mathematics, Graduate School of Engineering
\\
Hiroshima University, 
Higashihiroshima, 739-8527, Japan
}

\maketitle
\footnote[0]{E-mail: tshibata@hiroshima-u.ac.jp}
\footnote[0]{This work was supported by JSPS KAKENHI Grant Number JP17K05330.}

\begin{abstract}
We consider the nonlinear eigenvalue problem
$$
[D(u(t))u(t)']' + \lambda g(u(t)) = 0, \enskip u(t) > 0, 
\enskip t \in I := (0,1), \enskip u(0) = u(1) = 0,
$$
which comes from the porous media type 
equation. 
Here, $D(u) = pu^{2n} + \sin u$ 
($n \in \mathbb{N}$, $p > 0$: given constants), 
$g(u) = u$ or $g(u) = u + \sin u$. 
$\lambda > 0$ is a bifurcation parameter which is a continuous function of 
$\alpha = \Vert u_\lambda\Vert_\infty$ of the solution 
$u_\lambda$ 
corresponding to $\lambda$, and 
is expressed as $\lambda = \lambda(\alpha)$. 
Since our equation contains oscillatory term 
in diffusion term, it seems significant to study how 
this oscillatory term gives effect to 
the structure of bifurcation curves $\lambda(\alpha)$. 
We prove that the simplest case $D(u) = u^{2n} + \sin u$ and 
$g(u) = u$ gives us the most significant phenomena to the 
global behavior of $\lambda(\alpha)$. 

\vspace{0.5cm}

\noindent{MSC}: Primary 34C23, Secondary 34F10

\noindent
{\it Keywords}: Precise structure of bifurcation curves; Oscillatory 
nonlinear diffusion

\end{abstract}

\section{Introduction} 		      

We study the 
following nonlinear eigenvalue problems 
\begin{eqnarray}
[D(u(t))u(t)']' + \lambda g(u(t)) &=& 0, \quad t \in I:= (0,1),
\\
u(t) &>& 0, \quad t \in I,
\\
u(0) &=& u(1) = 0,
\end{eqnarray}
where $D(u):= pu^{2n} + \sin u$ 
($n \in \mathbb{N}$, $p > 0$: given constants), $g(u) = u$ or 
$g(u) = u + \sin u$, and $\lambda > 0$ is a bifurcation parameter. 
We assume the following condition (A.1).

\vspace{0.2cm}

\noindent
(A.1) $D(u) > 0$ for $u > 0$. 

\vspace{0.2cm}

\noindent{}
Under the condition 
(A.1), we know from [9] that for a given $\alpha > 0$, there 
is a unique slution pair $(u_\alpha, \lambda)$ of 
(1.1)--(1.3) satisfying $\alpha = \Vert u_\alpha\Vert_\infty$. 
Moreover, $\lambda$ is parameterized by $\alpha > 0$ as $\lambda(\alpha)$ and 
is continuous for $\alpha > 0$. 

The purpose 
of this paper is to show how 
the oscillatory diffusion term $D(u)$ gives effect to 
the structure of bifurcation curves $\lambda(\alpha)$. To clarify 
our intention, let $n = p = 1$ in (1.1) for simplicity. 
Then (A.1) is satisfied and we have 
the equation 
\begin{eqnarray}
[\{u(t)^2 + \sin u(t)\}u'(t)]' + \lambda u(t) &=& 0, \quad t \in I.
\end{eqnarray}
Another equations similar to (1.4) are
\begin{eqnarray}
[u(t)^2u(t)']' + \lambda (u(t) + \sin u(t)) &=& 0, \quad t \in I,
\end{eqnarray}
\begin{eqnarray}
[\{u(t)^2 + \sin u(t)\}u(t)']' + \lambda (u(t) + \sin u(t)) &=& 0, 
\quad t \in I.
\end{eqnarray}
{\bf Question A.} {\it Consider (1.j) with (1.2)--(1.3) ($j = 4,5,6$). 
Then can we distinguish 
the global structure of $\lambda(\alpha)$ for (1.4), (1.5) and (1.6) 
or not ?}

We explain the motivation and back ground of Question A. 
Bifurcation problems with $D(u) \equiv 1$ are one of the main interest 
in the study of differential equations, and many 
results have been established concerning 
the asymptotic behavior of bifurcation curves from 
mathematical point of view. 
We refer to [1--3,7,8,10,11,12] and the references therein. 
Besides, the bifurcation problems with nonlinear diffusion appear 
in the various fields. 
The case 
$D(u) = u^k$ ($k > 0$) appears as the porous media equation in 
material science and logistic type model equation in population dynamics. 
In the latter case, it implies 
that the diffusion rate $D(u)$ 
depends on both the population density $u$ and a parameter $1/\lambda$. 
We refer to [9, 13, 20] and the references therein. 
Added to these, there are several papers studying  
the asymptotic behavior of oscillatory bifurcation curves. 
We refer to [4,5,6,7,14--18] and the references therein. 

Recently, the following equation has been considered in [19].
\begin{eqnarray}
[D(u(t))u(t)']' + \lambda g(u(t)) &=& 0, \quad t \in I
\end{eqnarray}
with (1.2)--(1.3). 
Here, $D(u) = u^k$, $g(u) = u^{2m-k-1} + \sin u$, and 
$m \in \mathbb{N}$, $k$ ($0 \le k < 2m-1$) are given 
constants. In particular, if we put $m = k = 2$, then we have 
the equation (1.5). 
In [19], the following result has been obtained.

\vspace{0.2cm}

\noindent
{\bf Theorem 1.1 ([19]).} {\it Consider (1.7) with (1.2)--(1.3). Then 
as $\alpha \to \infty$,
\begin{eqnarray}
\lambda(\alpha) = 4m\alpha^{2k+2-2m}
\left\{A_{k,m}^2 - 2A_{k,m}\sqrt{\frac{\pi}{2m}}\alpha^{k+(1/2)-2m}
\sin\left(\alpha-\frac{\pi}{4}\right) + o(\alpha^{k+(1/2)-2m})\right\},
\end{eqnarray}
where 
\begin{eqnarray}
A_{k,m} = \int_0^1 \frac{s^k}{\sqrt{1-s^{2m}}}ds.
\end{eqnarray}
}

\vspace{0.1cm}

\noindent
If $m = k = 2$ in (1.8), then 
the asymptotic formula for $\lambda(\alpha)$ as $\alpha \to \infty$ 
for (1.5) is given by 
\begin{eqnarray}
\lambda(\alpha) &=& 8\alpha^2\left\{A_{2,2}^2
- A_{2,2}\sqrt{\pi}\alpha^{-3/2}\sin\left(\alpha - \frac{\pi}{4}\right) 
+ o(\alpha^{-3/2})\right\}.
\end{eqnarray}
Moreover, it was shown in [19] that if $m = k = 2$, then $\lambda(\alpha) 
= 4B_0^2\alpha^2(1 + o(1))$ as $\alpha \to 0$, where $B_0$ is a positive 
constant. We understand from Theorems 1.1 that, by the 
effect of $\sin u$, $\lambda(\alpha)$ oscillates and crosses the 
curve $\lambda = 8A_{2,2}\alpha^2$,
which are the {\it original} 
bifurcation curve obtained from the equation (1.5) 
{\it without} $\sin u$, 
infinitely many times. Motivated by this, 
we would like to study how oscillatory 
nonlinear diffusion 
influences on the structure of $\lambda(\alpha)$. 
The question we have to ask here is whether the 
oscillatory term $\sin u$ in $D(u)$ gives the same influence 
on the asymptotic behavior of $\lambda(\alpha)$ as (1.8) or not.

\noindent
Now we state our main results which give us 
the answer to this question. 

\vspace{0.2cm}

\noindent
{\bf Theorem 1.2.} {\it Assume (A.1). Consider (1.1)--(1.3) 
with $g(u) = u$.  

\noindent
{(i)} As $\alpha \to \infty$,
\begin{eqnarray}
\lambda(\alpha) &=& 
\frac{4(n+1)}{p}
\left(
p^2A_{2n,n+1}^2\alpha^{2n} + 2A_{2n,n+1}(p-1)\sqrt{\frac{\pi}{2(n+1)}}\alpha^{-1/2}
\sin\left(\alpha - \frac{\pi}{4}\right) 
\right.
\nonumber
\\
&&\left.\qquad \qquad \quad
+ O(\alpha^{-1})
\right).
\end{eqnarray}
In particular, consider (1.4) with (1.2)-(1.3). 
Then as $\alpha \to \infty$,
\begin{eqnarray}
\lambda(\alpha) = 8A_{2n,n+1}^2\alpha^2 + O(\alpha^{-1}).
\end{eqnarray}
(ii) Let $n = 1$. Then as $\alpha \to 0$,
\begin{eqnarray}
\lambda(\alpha) = 6\alpha(C_0^2 + 2pC_0C_{1}\alpha + O(\alpha^2)),
\end{eqnarray}
where
\begin{eqnarray}
C_0:= \int_0^1 \frac{s}{\sqrt{1-s^3}}ds, \quad 
C_{1}:= \int_0^1 \frac{1}{\sqrt{1-s^3}}\left(s^2 - 
\frac{3s(1-s^4)}{8(1-s^3)}\right)ds.
\end{eqnarray}
}

\vspace{0.2cm}

\noindent
Therefore, our conclusion is that 
the decay rate of the second term of (1.12) is different 
from (1.10). Namely, the answer to Question A is affirmative. The global structures of the bifurcation curves for 
(1.4) and (1.5) do not coincide each other. Certainly, 
the future direction of this study will be to obtain the 
exact second term of 
(1.11) when $p = 1$, although it seems very difficult to 
get it. 

\vspace{0.2cm}

\noindent
{\bf Remark.} (i) If $n = 1$ and $p \not=1$ in (1.11), 
then the rough image of the shape of $\lambda(\alpha)$ is like Fig. 1. 
However, $\lambda(\alpha)$ in (1.11) 
is closer to 
$\lambda = 8A_{2,2}^2\alpha^2$ than $\lambda(\alpha)$ in (1.8). 
The reason is as follows. The combination of the power 
nonlinearity $u^2$ in $D(u)$ and $\sin u$ in $g(u)$ in (1.5) 
gives stronger effect to 
the second term of $\lambda(\alpha)$ than that given by the 
combination of $\sin u$ in $D(u)$ 
and $u$ in $g(u)$ in (1.4).

\noindent
{(ii)} Theorem 1.2 (ii) is only proved for the case $n = 1$ to 
show that the rough picture of $\lambda(\alpha)$ is almost the same as 
that of Fig. 1 if $p \not=1$. 
Certainly, we easily obtain Theorem 1.2 (ii) for 
the case $n \ge 2$.

\vspace{0.2cm}

The following Theorem 1.3 gives us the negative answer to Question A.

\vspace{0.2cm}

\noindent
{\bf Theorem 1.3.} {\it 
Consider (1.6) with (1.2)--(1.3). 

\noindent
(i) The asymptotic formula (1.10) holds as $\alpha \to \infty$. 

\noindent
(ii) The following asymptotic formula holds as $\alpha \to 0$.
\begin{eqnarray}
\lambda(\alpha) = 3\alpha
\left(C_0^2 + 2C_0C_{1}\alpha + O(\alpha^2)\right). 
\end{eqnarray}
}

\vspace{0.2cm}

\noindent
We find from Theorems 1.2 and 1.3 that 
$\sin u$ in diffusion term 
has deep influences on the global behavior of $\lambda(\alpha)$. 

We prove Theorems 1.2 and 1.3 by using time-map method and 
stationary phase method. 

\section{Proof of Theorem 1.2 (i)}

In this section, let $D(u) = pu^{2n} + \sin u$, $g(u) = u$ and  
$\alpha \gg 1$. We denote by $C$ the various positive constants 
independent of $\alpha \gg 1$. We put 
\begin{eqnarray}
\Lambda:= \left\{\alpha > 0 \enskip \vert \enskip 
g(\alpha) > 0, \displaystyle{\int_u^\alpha} 
g(t)D(t) dt > 0 \enskip \mbox{for all} \enskip u \in [0,\alpha)\right\}.
\end{eqnarray}
It follows from [9, (2.7)] that if $\alpha \in \Lambda$, 
then $\lambda(\alpha)$ is well defined. 
By (A.1), we have 
$D(t) > 0$, $g(t) > 0$ for $t > 0$. So $g(t)D(t) > 0$ 
for $t > 0$ holds. 
Hence, $\Lambda \equiv \mathbb{R}_+$. 
By this and the generalized time-map in [9, (2.5)] 
(cf. (2.7) below) and the 
time-map argument in [8, Theorem 2.1], we find that 
for any given $\alpha > 0$, 
there is a unique solution pair 
$(u_\alpha, \lambda) \in C^2(I) \bigcap C(\bar{I}) \times \mathbb{R}_+$ 
of (1.1)--(1.3) satisfying 
$\alpha= \Vert u_\alpha\Vert_\infty$. Moreover, 
$\lambda$ is parameterized by $\alpha$ as 
$\lambda = \lambda(\alpha)$ and is a 
continuous function for $\alpha > 0$. 
It is well known that if 
$(u_\alpha, \lambda(\alpha)) \in C^2({I}) \bigcap C(\bar{I}) 
\times \mathbb{R}_+$ satisfies 
(1.1)--(1.3), then 
\begin{eqnarray}
&&u_\alpha(t) = u_\alpha(1-t), \quad 0 \le t \le 1, \\
&&u_\alpha\left(\frac12\right) = \max_{0 \le t \le 1}u_\alpha(t) 
= \alpha, \\
&&u_\alpha'(t) > 0, \quad 0 < t < \frac12.
\end{eqnarray}
We put 
\begin{eqnarray}
G(u)&:=& \int_0^u D(x)g(x)dx = \int_0^u (px^{2n} + \sin x)xdx 
\\
&=& \frac{p}{2n+2}u^{2n+2} - u\cos u + \sin u,
\nonumber
\\
M_1&:=& \alpha\cos\alpha - \alpha s\cos(\alpha s), 
\quad M_2:= \sin\alpha - \sin(\alpha s) \quad (0 \le s \le 1).
\end{eqnarray}
For $0 \le s \le 1$ and $\alpha \gg 1$, we have 
\begin{eqnarray}
\frac{\vert M_1 \vert 
+ \vert M_2 \vert}{\alpha^{2n+2}(1-s^{2n+2})} 
\le C\alpha^{-2n}.
\end{eqnarray}
By this, Taylor expansion and putting $u = s\alpha$, we have from 
[9] that 
\begin{eqnarray}
\sqrt{\frac{\lambda}{2}} 
&=& \int_0^\alpha \frac{D(u)}
{\sqrt{G(\alpha)-G(u)}}du
\\
&=& \alpha
\int_0^1 \frac{p\alpha^{2n}s^{2n} + \sin(\alpha s)}
{
\sqrt{\frac{p}{2n+2}\alpha^{2n+2}(1-s^{2n+2})
-(M_1 - M_2)}
}ds
\nonumber
\\
&=& \sqrt{\frac{2n+2}{p}}\alpha^{-n}
\int_0^1 \frac{p\alpha^{2n}s^{2n} + \sin(\alpha s)}
{\sqrt{1-s^{2n+2}}
\sqrt{1 - \frac{2n+2}{p\alpha^{2n+2}(1-s^{2n+2})}
\{M_1 - M_2\}
}
}ds
\nonumber
\\
&=& 
\sqrt{\frac{2n+2}{p}}\alpha^{-n}
\int_0^1 \frac{p\alpha^{2n}s^{2n} + \sin(\alpha s)}
{\sqrt{1-s^{2n+2}}}
\nonumber
\\
&&\qquad \qquad \qquad \times
\left(1 + \frac{n+1}{p\alpha^{2n+2}(1-s^{2n+2})}
\{M_1 - M_2\}(1 + O(\alpha^{-2n}))
\right)ds.
\nonumber
\end{eqnarray}
This implies that 
\begin{eqnarray}
\sqrt{\frac{\lambda}{2}} &=& 
\sqrt{\frac{2n+2}{p}}\alpha^{-n}\{J_1 + J_2 + J_3 + J_4 + J_5 
\}(1+ O(\alpha^{-2n})),
\end{eqnarray}
where
\begin{eqnarray}
J_1&:=& p\alpha^{2n}\int_0^1 \frac{s^{2n}}{\sqrt{1-s^{2n+2}}}ds 
= pA_{2n,n}\alpha^{2n}, 
\\
J_2&:=& \int_0^1 \frac{\sin(\alpha s)}{\sqrt{1-s^{2n+2}}}ds,
\\
J_3&:=& \frac{n+1}{p\alpha}\int_0^1 \frac{s^{2n}}{(1-s^{2n+2})^{3/2}}
(\cos\alpha - s\cos(\alpha s))ds,
\\
J_4&:=&  -\frac{n+1}{p\alpha^2}\int_0^1 \frac{s^{2n}}{(1-s^{2n+2})^{3/2}}
(\sin\alpha - \sin(\alpha s)ds,
\\
J_5 &:=& \frac{n+1}{p\alpha^{2n+2}}
\int_0^1 \frac{\sin(\alpha s)}{(1-s^{2n+2})^{3/2}}
\{M_1-M_2\}ds.
\end{eqnarray}

\vspace{0.2cm}

\noindent
To calculate $J_2 \sim J_5$, we use the following equality.

\vspace{0.2cm}

\noindent
{\bf Lemma 2.1 ([5, Lemma 2], [7, Lemma 2.25]).} 
{\it Assume that the function 
$f(r) \in C^2[0,1]$, $w(r) = \cos(\pi r/2)$. 
Then as $\mu \to \infty$
\begin{eqnarray}
\int_0^1 f(r)e^{i\mu w(r)}dr = e^{i(\mu -(\pi/4))}
\sqrt{\frac{2}{\mu\pi}}f(0) + O\left(\frac{1}{\mu}\right).
\end{eqnarray}
In particular, by taking the real and 
imaginary parts of (2.15), as $\mu \to \infty$, 
\begin{eqnarray}
\int_0^1 f(r)\cos(\mu w(r))dr &=& 
\sqrt{\frac{2}{\mu\pi}}f(0)
\cos\left(\mu-\frac{\pi}{4}\right) 
+ O\left(\frac{1}{\mu}\right),
\\
\int_0^1 f(r)\sin(\mu w(r))dr &=& 
\sqrt{\frac{2}{\mu\pi}}f(0)
\sin\left(\mu-\frac{\pi}{4}\right) 
+ O\left(\frac{1}{\mu}\right).
\end{eqnarray}
}

\vspace{0.2cm}

\noindent
{\bf Lemma 2.2.} {\it As $\alpha \to \infty$,}
\begin{eqnarray}
J_2 =\sqrt{\frac{\pi}{2(n+1)\alpha}
}\sin\left(\alpha - \frac{\pi}{4}\right) + O(\alpha^{-1}). 
\end{eqnarray}
{\it Proof.} Putting $s = \sin\theta$, 
$\theta = \frac{\pi}{2}(1-x)$ and using Lemma 2.1, we have
\begin{eqnarray}
J_2 &=& \int_0^1 \frac{\sin(\alpha s)}
{\sqrt{1-s^2}\sqrt{1 + s^2 + \cdots + s^{2n}}}ds
\\
&=& \int_0^{\pi/2} \frac{1}{\sqrt{1 + \sin^2\theta + \cdots + \sin^{2n}\theta}}
\sin(\alpha\sin\theta)d\theta
\nonumber
\\
&=& \frac{\pi}{2}\int_0^{\pi/2}
\frac{1}
{
\sqrt{1 + 
\cos^{2}\left(\frac{\pi}{2}x\right) + \cdots 
+ \cos^{2n}\left(\frac{\pi}{2}x\right)}
}
\sin\left(
\alpha\cos\left(\frac{\pi}{2}x
\right)\right)dx
\nonumber
\\
&=& \sqrt{\frac{\pi}{2(n+1)\alpha}
}\sin\left(\alpha - \frac{\pi}{4}\right) + O(\alpha^{-1}).
\nonumber
\end{eqnarray}
Thus the proof is complete. \qed

\vspace{0.2cm}

\noindent
{\bf Lemma 2.3.} {\it As $\alpha \to \infty$,}
\begin{eqnarray}
J_3 = -\frac{1}{p}\sqrt{\frac{\pi}{2(n+1)\alpha}}
\left\{\sin\left(\alpha-\frac{\pi}{4}\right) 
- \alpha^{-1}
\cos\left(\alpha - \frac{\pi}{4}\right)\right\} + O(\alpha^{-1}).
\end{eqnarray}
{\it Proof.} We put $J_3 = (n+1)J_{31}/(p\alpha)$, $s = \sin\theta$ and 
$K(\theta) := \sin^{2n}\theta/(1 + \sin^2\theta + \cdots 
+ \sin^{2n}\theta)^{3/2}$. 
Then by integration by parts, 
\begin{eqnarray}
J_{31} &=& \int_0^1 \frac{s^{2n}}{(1-s^2)^{3/2}
(1 + s^2 + \cdots + s^{2n})^{3/2}}(\cos \alpha - s\cos(\alpha s))ds
\\
&=& \int_0^{\pi/2}\frac{1}{\cos^2\theta}K(\theta)
(\cos\alpha - \sin\theta\cos(\alpha\sin\theta))d\theta
\nonumber
\\
&=& \left[\tan\theta K(\theta)(\cos\alpha - \sin\theta\cos(\alpha\sin\theta))
\right]_0^{\pi/2} 
\nonumber
\\
&&\mbox{}-\int_0^{\pi/2} \tan\theta K'(\theta)
(\cos\alpha - \sin\theta\cos(\alpha\sin\theta))d\theta 
\nonumber
\\
&&\mbox{}- \int_0^{\pi/2}\tan\theta K(\theta)
\{-\cos\theta\cos(\alpha\sin\theta)+\alpha\sin\theta\cos\theta
\sin(\alpha\sin\theta)\}d\theta
\nonumber
\\
&:=& J_{311} - J_{312} + J_{313}.
\nonumber
\end{eqnarray} 
By using l'H\^opital's rule, we have 
\begin{eqnarray}
&&\lim_{\theta \to \pi/2}\frac{\cos\alpha - \sin\theta\cos(\alpha\sin\theta)}
{\cos\theta} 
\\
= &&\lim_{\theta \to \pi/2}\frac{-\cos\theta\cos(\alpha\sin\theta) 
+ \alpha\sin\theta\cos\theta
\sin(\alpha\sin\theta)}{-\sin\theta} = 0.
\nonumber
\end{eqnarray}
We see from this that $J_{311} = 0$. 
Moreover, by direct calculation, we see that $J_{312} = O(1)$. Now, putting 
$\theta = \frac{\pi}{2}(1-x)$ and using Lemma 2.1, we have
\begin{eqnarray}
J_{313} &=& \int_0^{\pi/2}\sin\theta K(\theta)\cos(\alpha\sin\theta)d\theta 
-\alpha\int_0^{\pi/2}K(\theta)\sin^2\theta\sin(\alpha\sin\theta)d\theta
\\
&=& \frac{\pi}{2}\int_0^{1} \frac{\cos^{2n+1}\left(\frac{\pi}{2}x\right)}
{\left(1 + \cos^2\left(\frac{\pi}{2}x\right) + \cdots 
+ \cos^{2n}\left(\frac{\pi}{2}x\right)\right)^{3/2}}
\cos\left(\alpha\cos\left(\frac{\pi}{2}x\right)\right)dx
\nonumber
\\
&&\mbox{}
-\frac{\pi}{2}\alpha\int_0^{1} \frac{\cos^{2n+2}\left(\frac{\pi}{2}x\right)}
{\left(1 + \cos^2\left(\frac{\pi}{2}x\right) + \cdots 
+ \cos^{2n}\left(\frac{\pi}{2}x\right)\right)^{3/2}}
\sin\left(\alpha\cos\left(\frac{\pi}{2}x\right)\right)dx
\nonumber
\\
&=& 
(n+1)^{-3/2}\sqrt{\frac{\pi}{2\alpha}}
\cos\left(\alpha - \frac{\pi}{4}\right) 
- (n+1)^{-3/2}\sqrt{\frac{\pi\alpha}{2}}
\sin\left(\alpha - \frac{\pi}{4}\right) 
+ O(1).
\nonumber
\end{eqnarray}
This implies (2.20). Thus the proof is complete. \qed

\vspace{0.2cm}

\noindent
{\bf Lemma 2.4.} {\it As $\alpha \to \infty$,}
\begin{eqnarray}
J_4 = -\frac{1}{p}\sqrt{\frac{\pi}{2(n+1)\alpha}}\alpha^{-1}
\cos\left(\alpha - \frac{\pi}{4}
\right) + O(\alpha^{-2}).
\end{eqnarray}
{\it Proof.} We put $J_4 = -(n+1)J_{41}/(p\alpha^2)$ and $s = \sin\theta$. 
Then by the same argument as that to obtain $J_{31}$ in (2.23), 
we obtain 
\begin{eqnarray}
J_{41} &=& \int_0^1 \frac{s^{2n}}{(1-s^{2n+2})^{3/2}}
(\sin\alpha - \sin(\alpha s)ds
\\
&=& \int_0^{\pi/2} \frac{1}{\cos^2\theta}K(\theta)
(\sin\alpha - \sin(\alpha\sin\theta))d\theta
\nonumber
\\
&=& \left[\tan\theta K(\theta)(\sin\alpha - \sin(\alpha\sin\theta))
\right]_0^{\pi/2} 
\nonumber
\\
&&\mbox{}- \int_0^{\pi/2} \tan\theta K'(\theta)
(\sin\alpha - \sin(\alpha\sin\theta))d\theta
\nonumber
\\
&&\mbox{}+\alpha\int_0^{\pi/2}\sin\theta K(\theta)\cos(\alpha\sin\theta)d\theta
\nonumber
\\
&=& \frac{\pi}{2}\alpha \int_0^{1} 
\frac{\cos^{2n+1}\left(\frac{\pi}{2}x\right)}
{\left(1 + \cos^2\left(\frac{\pi}{2}x\right) + \cdots + 
\cos^{2n}\left(\frac{\pi}{2}x\right)\right)^{3/2}}
\cos\left(\alpha\cos\left(\frac{\pi}{2}x\right)\right)dx + O(1)
\nonumber
\\
&=& (n+1)^{-3/2}\sqrt{\frac{\pi\alpha}{2}}\cos\left(\alpha - \frac{\pi}{4}
\right) + O(1).
\nonumber
\end{eqnarray}
By this, we obtain (2.24). Thus the proof is complete. \qed

\vspace{0.2cm}

\noindent
{\it Proof of Theorem 1.2 (i)}. By (2.7) and (2.14), we see that 
$J_5 = O(\alpha^{-2n})$. 
By this, (2.9) and Lemmas 2.2--2.4, 
we obtain 
\begin{eqnarray}
\sqrt{\frac{\lambda}{2}} = \sqrt{\frac{2n+2}{p}}
\alpha^{-n}
\left[pA_{2n,n}\alpha^{2n} + 
\left(1 - \frac{1}{p}\right)\sqrt{\frac{\pi}{2(n+1)\alpha}}
\sin\left(\alpha - \frac{\pi}{4}\right) + O(\alpha^{-1})\right].
\end{eqnarray}
By this, we obtain 
\begin{eqnarray}
\lambda 
&=& \frac{4(n+1)}{p}\alpha^{-2n}
\left[p^2A_{2n,n+1}^2\alpha^{4n} + 
2\left(p - 1\right)A_{2n,n+1}\alpha^{2n-(1/2)}\sqrt{\frac{\pi}{2(n+1)}}
\sin\left(\alpha - \frac{\pi}{4}\right) 
\right.
\nonumber
\\
&&\qquad \left.\qquad \qquad \qquad + O(\alpha^{2n-1})\right].
\end{eqnarray}
This implies Theorem 1.2 (i). Thus the proof is complete. \qed 

\section{Proof of Theorem 1.3 (i)}

In this section, let $D(u) = u^2 + \sin u$ and $g(u) = u + \sin u$. 
It follows from [9, (2.7)], we also find, as in Section 2, that 
for any given $\alpha > 0$, 
there is a unique classical solution pair 
$(\lambda, u_\alpha)$ of (1.1)--(1.3) satisfying 
$\alpha= \Vert u_\alpha\Vert_\infty$. Moreover, 
$\lambda$ is parameterized by $\alpha$ as 
$\lambda = \lambda(\alpha)$ and is a 
continuous function for $\alpha > 0$. 
Let $u \ge 0$. We put 
\begin{eqnarray}
G(u) &:=& \int_0^u g(y)D(y)dy 
\\
&=& \frac{1}{4}u^{4} 
- u \cos u + \sin u + (2u\sin u - (u^2-2)\cos u - 2) + 
\frac12\left(u - \frac12\sin 2u\right) 
\nonumber
\\
&:=& \frac{1}{4}u^4 + G_1(u).
\nonumber
\end{eqnarray}
For $0 \le s \le 1$ and $\alpha \gg 1$, we have 
\begin{eqnarray}
G(\alpha) - G(\alpha s) &=& \frac{1}{4}\alpha^4(1-s^4) 
+ G_1(\alpha) - G_1(\alpha s) 
\\
&=& \frac{1}{4}\alpha^4(1-s^4)
- (\alpha\cos\alpha - \alpha s\cos(\alpha s)) 
+ (\sin \alpha - \sin (\alpha s)) 
\nonumber
\\
&&\mbox{} + 2(\alpha\sin\alpha - \alpha s\sin(\alpha s)) 
- (\alpha^2\cos\alpha -\alpha^2 s^2\cos(\alpha s)) 
\nonumber
\\
&&\mbox{} + 2(\cos\alpha - \cos(\alpha s)) + \frac12\alpha(1-s) 
- \frac14(\sin2\alpha - \sin2\alpha s))
\nonumber
\\
&:=& \frac{1}{4}\alpha^4(1-s^4) 
- I_1 + I_2 + I_3 - I_4 + I_5 + I_6 - I_7. 
\nonumber
\end{eqnarray}
It is easy to see that for $0 \le s \le 1$, 
\begin{eqnarray}
&&\left\vert \frac{I_4}
{\alpha^{4}(1-s^{4})}\right\vert \le C\alpha^{-1}, 
\\
&&\left\vert \frac{I_1}
{\alpha^{4}(1-s^{4})}\right\vert, 
\left\vert \frac{I_3}
{\alpha^{4}(1-s^{4})}\right\vert \le C\alpha^{-2},
\\
&&
\left\vert \frac{I_2}
{\alpha^{4}(1-s^{4})}\right\vert, 
\left\vert \frac{I_5}
{\alpha^{4}(1-s^{4})}\right\vert, 
\left\vert \frac{I_6}
{\alpha^{4}(1-s^{4})}\right\vert,  
\left\vert \frac{I_7}
{\alpha^{4}(1-s^{4})}\right\vert \le C\alpha^{-3}. 
\end{eqnarray}
By putting $u = \alpha s$, 
(3.3)--(3.5) and Taylor expansion, we have from [9, (2.5)] that 
\begin{eqnarray}
\sqrt{\frac{\lambda(\alpha)}{2}} &=& \int_0^\alpha \frac{D(u)}
{\sqrt{G(\alpha)-G(u)}}du
\\
&=& \int_0^\alpha \frac{u^2 + \sin u}{\sqrt{\frac{1}{4}(\alpha^{4}-u^{4}) 
+ G_1(\alpha)-G_1(u)}}du
\nonumber
\\
&=& \alpha\int_0^1 \frac{\alpha^2s^2 + \sin\alpha s}
{\sqrt{\alpha^4(1-s^{4})/4 + G_1(\alpha)-G_1(\alpha s)}}ds
\nonumber
\\
&=& 
2\alpha^{-1}
\int_0^1 \frac{\alpha^2s^2 + \sin\alpha s}
{\sqrt{1-s^{4}}\sqrt{1 + \frac{4}{\alpha^4(1-s^4)}(G_1(\alpha)-G_1(\alpha s))}}ds
\nonumber
\\
&=& 2\alpha^{-1}\int_0^1 
\frac{1}{\sqrt{1-s^4}}(\alpha^2 s^2 + \sin\alpha s)
\nonumber
\\
&&\times
\left\{1 - \frac{2}{\alpha^4(1-s^4)}(G_1(\alpha) - G_1(\alpha s))
(1 + O(\alpha^{-1}))\right 
\}ds.
\nonumber
\end{eqnarray}

\vspace{0.2cm}

\noindent
Now we show that the leading and second terms of 
the right hand side of (3.6) are 
\begin{eqnarray}
L_1&:=& 2\alpha\int_0^1 \frac{s^2}{\sqrt{1-s^4}}ds 
= 2A_{2,2}\alpha, 
\\
L_4&:=& 4\alpha^{-5}\int_0^1 \frac{s^2}{(1-s^4)^{3/2}}I_4 ds.
\end{eqnarray}
Indeed, by (3.3)--(3.5), we obtain 
\begin{eqnarray}
2\alpha^{-1}\int_0^1 \frac{\vert \alpha^2s^2 + \sin\alpha s\vert}
{\sqrt{1-s^4}}\cdot\frac{1}{\alpha^4(1-s^4)}
\vert I_1 + I_3 + I_5 + I_6 + I_7\vert ds 
= O(\alpha^{-1}).
\end{eqnarray}
Furthermore, by Lemma 2.2, 
\begin{eqnarray}
L_2&:=& 2\alpha^{-1}\int_0^1 \frac{\sin(\alpha s)}
{\sqrt{1-s^4}}ds = (1 + o(1))\sqrt{\pi}\alpha^{-3/2}
\sin\left(\alpha - \frac{\pi}{4}\right).
\end{eqnarray}
\vspace{0.2cm}

\noindent
We calculate $L_4$ by Lemma 2.1. 

\vspace{0.2cm}

\noindent
{\bf Lemma 3.1.} {\it As $\alpha \to \infty$,}
\begin{eqnarray}
L_4 = 
- \sqrt{\pi}\alpha^{-1/2}\sin\left(\alpha - \frac{\pi}{4}\right)
+ O(\alpha^{-1}).
\end{eqnarray}
{\it Proof.} We put $s = \sin\theta$. Then 
\begin{eqnarray}
L_4 &=& 4\alpha^{-1}\int_0^1\frac{s^2(\cos\alpha - s^2\cos(\alpha s))}
{(1-s^2)^{3/2}(1 + s^2)^{3/2}}ds
\\
&=& 4\alpha^{-1}\int_0^{\pi/2} \frac{1}{\cos^2\theta}
Y(\theta)(\cos\alpha - \sin^2\theta\cos(\alpha\sin\theta))
d\theta,
\nonumber
\end{eqnarray}
where $Y(\theta) = \sin^2\theta /
(1 + \sin^2\theta)^{3/2}$. By Integration by parts, we have 
\begin{eqnarray}
L_4 &=& 4\alpha^{-1}\left[\tan\theta Y(\theta)
(\cos\alpha - \sin^2\theta\cos(\alpha\sin\theta))\right]_0^{\pi/2}
\\
&&\mbox{}
- 4\alpha^{-1}\int_0^1 \tan\theta\left\{Y(\theta)
(\cos\alpha - \sin^2\theta\cos(\alpha\sin\theta))\right\}'d\theta 
\nonumber
\\
&=& 4\alpha^{-1}(L_{41} - L_{42}).
\nonumber
\end{eqnarray}
By using l'H\^opital's rule, we have 
\begin{eqnarray}
&&\lim_{\theta \to \pi/2}
\frac{\cos\alpha - \sin^2\theta\cos(\alpha\sin\theta)}{\cos\theta}
\\
=
&&\lim_{\theta \to \pi/2}\frac{-2\sin\theta\cos\theta
\cos(\alpha\sin\theta) + \alpha\sin^2\theta
\cos\theta\sin(\alpha\sin\theta)}{-\sin\theta} = 0.
\nonumber
\end{eqnarray}
By this, we see that $L_{41} = 0$. It is easy to see that 
\begin{eqnarray}
\int_0^1 \tan\theta\left\{Y(\theta)'
(\cos\alpha - \sin^2\theta\cos(\alpha\sin\theta))\right\}d\theta = O(1).
\end{eqnarray}
By this and putting $\theta = \frac{\pi}{2}(1-x)$ and using Lemma 2.1, we have 
\begin{eqnarray}
L_{42} &=& \int_0^1 \tan\theta Y(\theta)
(\cos\alpha - \sin^2\theta\cos(\alpha\sin\theta))'d\theta + O(1)
\\
&=& \alpha\int_0^{\pi/2} \frac{\sin^5\theta}{(1+\sin^2\theta)^{3/2}}
\sin(\alpha\sin\theta)d\theta + O(1)
\nonumber
\\
&=&
\frac{\pi}{2}\alpha\int_0^1 \frac{\cos^5\left(\frac{\pi}{2}x\right)}
{(1 + \cos^2\left(\frac{\pi}{2}x\right))^{3/2}}
\sin\left(\alpha\cos\left(\frac{\pi}{2}x\right)\right)dx + O(1)
\nonumber
\\
&=& \frac{\sqrt{\pi\alpha}}{4}\sin\left(\alpha - \frac{\pi}{4}\right)
+ O(1).
\nonumber
\end{eqnarray}
By this and (3.13), we obtain (3.11). Thus the proof is complete. \qed

\vspace{0.2cm}

\noindent
{\it Proof of Theorem 1.3 (i).} 
By (3.6), (3.7), Lemma 3.1, we obtain 
\begin{eqnarray}
\sqrt{\frac{\lambda}{2}} &=& 2A_{2,2}\alpha 
- \sqrt{\pi}\alpha^{-1/2}\sin\left(\alpha - \frac{\pi}{4}\right) 
+ o(\alpha^{-1/2}).
\end{eqnarray}
This implies (1.10). Thus the proof 
of Theorem 1.4 (i) is complete. \qed

\section{Proofs of Theorems 1.2 (ii) and 1.3 (ii)}
In this section, let $0 < \alpha \ll 1$.

\noindent
{\it Proof of Theorem 1.2 (ii).} Let $n = 1$, namely, 
$D(u) = pu^2 + \sin u$. By (2.5), Taylor expansion and 
direct calculation, 
for $0 \le s \le 1$, we have 
\begin{eqnarray}
G(\alpha) - G(\alpha s) &=& \frac13\alpha^3(1-s^3) 
+ \frac{p}{4}\alpha^4(1-s^4) + O(\alpha^5)
(1-s^5).
\end{eqnarray}
By this, outting $\theta= \alpha s$ and (2.8), we obtain 
\begin{eqnarray}
\sqrt{\frac{\lambda}{2}} &=& 
\alpha\int_0^1 \frac{\alpha s + p\alpha^2s^2 + O(\alpha^3)}
{\sqrt{\frac13\alpha^3(1-s^3) 
+ \frac{p}{4}\alpha^4(1-s^4) + O(\alpha^5)
(1-s^5)}}ds
\\
&=& \sqrt{3\alpha}\int_0^1 
\frac{s + p\alpha s^2 + O(\alpha^2)}
{\sqrt{(1-s^3) + \frac{3p}{4}\alpha(1-s^4) + O(\alpha^2)(1-s^5)}}ds
\nonumber
\\
&=& \sqrt{3\alpha}\int_0^1 \frac{1}{\sqrt{1-s^3}}
(s + p\alpha s^2 + O(\alpha^2))
\left(1 - \frac{3p(1-s^4)}{8(1-s^3)}\alpha + O(\alpha^2)\right)ds
\nonumber
\\
&=& \sqrt{3\alpha}\left\{\int_0^1 \frac{s}{\sqrt{1-s^3}}ds 
+ p\alpha\int_0^1 \frac{1}{\sqrt{1-s^3}}
\left(s^2-\frac{3s(1-s^4)}{8(1-s^3)}\right)ds 
+ O(\alpha^2)\right\}.
\nonumber
\end{eqnarray}
By this, we obtain Theorem 1.2 (ii). Thus the proof is complete. \qed

\noindent
{\it Proof of Theorem 1.3 (ii).} By (3.2), Taylor expansion and 
direct calculation, 
for $0 \le s \le 1$, we have 
\begin{eqnarray}
G(\alpha) - G(\alpha s) &=& \frac23\alpha^3(1-s^3) 
+ \frac12\alpha^4(1-s^4) + O(\alpha^5)
(1-s^5).
\end{eqnarray}
By this, putting $\theta = \alpha s$ and (3.2), we obtain 
\begin{eqnarray}
\sqrt{\frac{\lambda}{2}} 
&=& \alpha\int_0^1 
\frac{\alpha^2 s^2 + \sin(\alpha s)}
{\sqrt{\frac23\alpha^3(1-s^3) 
+ \frac12\alpha^4(1-s^4) + O(\alpha^5)(1-s^5)}}ds
\\
&=& \sqrt{\frac{3}{2\alpha}}\int_0^1 
\frac{\alpha s + \alpha^2s^2 + O(\alpha^3)s^3}
{\sqrt{1-s^3} 
\sqrt{1 + \frac34\alpha\frac{1-s^4}{1-s^3} + O(\alpha^2)
\frac{1-s^5}{1-s^3}}}ds
\\
\nonumber
\\
&=& \sqrt{\frac{3\alpha}{2}}\int_0^1 \frac{1}{\sqrt{1-s^3}}
(s + \alpha s^2 + O(\alpha^2))
\left(1 - \frac38\alpha\frac{1-s^4}{1-s^3} + O(\alpha^2)\right)ds
\nonumber
\\
&=& \sqrt{\frac{3\alpha}{2}}
\left(\int_0^1 \frac{s}{\sqrt{1-s^3}}ds 
+ \alpha
\left\{\int_0^1 \frac{s^2}{\sqrt{1-s^3}}ds 
-\frac38\int_0^1\frac{s(1-s^4)}{(1-s^3)^{3/2}}ds \right\} + O(\alpha^2)\right).
\nonumber
\end{eqnarray}
By this, we obtain Theorem 1.3 (ii). Thus the proof is complete. \qed

\end{document}